\documentclass[12pt]{article} 
\usepackage{amsmath, amssymb,latexsym, epsfig, color}
\input epsf

\newtheorem{theorem}{Theorem}[section]
\newtheorem{proposition}[theorem]{Proposition}

 \newtheorem{lemma}[theorem]{Lemma}

\def\proof{\smallskip\noindent {\it Proof: \ }} 
\def\endproof{\hfill$\square$\medskip}
 \def\Z{\mathbb{Z}} 
\def\C{\mathbb{C}} 
 \def\R{\mathbb{R}} 
 
\def\S{\mathbb{S}}

\def\B{\mathcal{B}}

\newcommand{\MON}{\mbox{\upshape \MON}}

\newcommand{\conv}{\mbox{\upshape conv}}

\newcommand{\fmax}{\mbox{\upshape fmax}}
\newcommand{\SM}{\mbox{\upshape SM}}
\newcommand{\Int}{\mbox{\upshape int}\,}
\newcommand{\vol}{\mbox{\upshape vol}\,}
\newcommand{\ctg}{\mbox{\upshape ctg}\,}

\title{A centrally symmetric version of the cyclic polytope}
\author{Alexander Barvinok
\thanks{Research partially supported by NSF grant DMS 0400617}\\
\small Department of Mathematics, \\[-0.8ex]
\small University of Michigan, Ann Arbor, Michigan 48109-1043, USA, \\[-0.8ex] 
\small \texttt{barvinok@umich.edu} 
\and
Isabella Novik 
\thanks{Research partially supported by Alfred P.~Sloan Research
Fellowship and NSF grant DMS-0500748 }\\ 
\small Department of Mathematics, Box 354350\\[-0.8ex] 
\small University of Washington, Seattle, WA 98195-4350, USA,\\[-0.8ex] 
\small \texttt{novik@math.washington.edu} }

\begin{document}

\maketitle

\begin{abstract}
We define a centrally symmetric analogue of the cyclic 
polytope and study its facial structure. 
 We conjecture that our polytopes provide asymptotically the largest 
 number of faces in all dimensions among all centrally symmetric polytopes 
with $n$ vertices of a given even dimension $d=2k$ when $d$ is fixed 
and $n$ grows. For a fixed even dimension $d=2k$
 and an integer $1 \leq j <k$ we prove that the maximum possible number 
of $j$-dimensional faces of a centrally symmetric  $d$-dimensional 
polytope with $n$ vertices is at least 
$\left(c_j(d)+o(1)\right) {n \choose j+1}$ for some $c_j(d)>0$ 
and at most $\left(1-2^{-d}+o(1)\right){n \choose j+1}$ as $n$ grows. 
We show that
 $c_1(d) \geq (d-2)/(d-1)$.\end{abstract}

\section{Introduction and main results}
To characterize the numbers that  arise as the face numbers of simplicial 
complexes of various types is a problem that has intrigued many researchers
over the last half century and has been solved
for quite a few classes of complexes, among them the class of all simplicial 
complexes \cite{Kr, Ka} as well as the class of all simplicial polytopes
\cite{BiLee, St80}. One of the precursors of the latter result
was the Upper Bound Theorem (UBT, for short)  \cite{McMu70} 
that provided  sharp upper bounds on the face numbers of {\em all}
$d$-dimensional polytopes with $n$ vertices. 
While the UBT is a classic by now, 
the situation for centrally symmetric polytopes  is wide open. 
For instance,  the largest number of edges, $\fmax(d,n;1)$, that a
$d$-dimensional centrally symmetric polytope on $n$ vertices can have
is unknown even for $d=4$ and no even conjectural bounds on this number exist.
 In this paper, we establish certain bounds 
on $\fmax(d,n; 1)$ and, more generally, on $\fmax(d,n;j)$, the maximum number 
of $j$-dimensional faces of a centrally symmetric $d$-dimensional polytope
 with $n$ vertices.
For every even dimension $d$ we construct a centrally symmetric polytope
with $n$ vertices, which, 
we conjecture, provides asymptotically the largest number of 
faces in every dimension
as $n$ grows and $d$ is fixed among all $d$-dimensional centrally 
symmetric polytopes with $n$ vertices.

Let us recall the basic definitions. A polytope will always mean
a convex polytope (that is, the convex hull of finitely many points), 
and a $d$-polytope---a $d$-dimensional polytope.
 A polytope $P\subset \R^d$ is 
{\em centrally symmetric} (cs, for short)
if for every $x\in P$, $-x$ belongs to $P$ as well, that is, $P=-P$. 
The number of $i$-dimensional faces ($i$-faces, for short) of $P$ is denoted
$f_i=f_i(P)$ and is called the {\em $i$th face number} of $P$.

The UBT proposed by Motzkin in 1957 \cite{Motz} and proved by
McMullen \cite{McMu70} asserts that among all 
$d$-polytopes  with $n$ vertices, 
the cyclic polytope, $C_d(n)$, maximizes the number of $i$-faces for
every $i$. Here the {\em cyclic polytope}, $C_d(n)$,
is the convex hull of $n$ distinct points on the 
{\em moment curve} $(t, t^2, \ldots, t^d)\in \R^d$ or on the
{\em trigonometric moment curve} 
$(\cos t, \sin t, \cos 2t, \sin 2t, \ldots, \cos kt, \sin kt)\in \R^{2k}$
(assuming $d=2k$). Both types of cyclic polytopes 
were investigated by Carath\'eodory \cite{Car} and later by Gale \cite{Gale}
who, in particular, showed that the two types are combinatorially 
equivalent  (for even $d$) and independent of the choice of points.
Cyclic polytopes were also rediscovered by Motzkin \cite{Motz, GrMotz}
and many others. We refer the readers to \cite{Barv} and \cite{Zieg} 
for more information on these amazing polytopes.

Here we define and study a natural centrally symmetric analog
of cyclic polytopes -- bicyclic polytopes.

\subsection{The symmetric moment curve and bicyclic polytopes}\label{curve}

Let us consider the curve 
$$
\SM_{2k}(t)= 
\Bigl(\cos t,\ \sin t,\ \cos 3t,\ \sin 3t,\ \ldots,
 \ \cos (2k-1)t,\ \sin (2k-1)t\Bigr) 
\quad \text{for} \quad t \in \R,$$
which we call the {\em symmetric moment curve}, 
$\SM_{2k}(t) \in \R^{2k}$. The difference between
$\SM_{2k}$ and the trigonometric moment curve 
is that we employ only odd multiples of $t$ in the former.
Clearly, $\SM_{2k}(t +2 \pi)=\SM_{2k}(t)$, so $\SM_{2k}$ defines a map
$\SM_{2k}: \R/2 \pi \Z \longrightarrow \R^{2k}$.
It is convenient to identify the quotient $\R/2 \pi \Z$ 
with the unit circle $\S^1 \subset \R^2$
via the map $t \longmapsto \left(\cos t, \sin t \right)$. 
In particular, $t$ and $t+\pi$ 
form a pair of antipodal points in $\S^1$. We observe that 
$$\SM_{2k}(t+\pi)=-\SM_{2k}(t),$$
so the symmetric moment curve $\SM_{2k}\left(\S^1\right)$
 is centrally symmetric about the origin.

Let $X\subset \S^1$ be a  finite set. 
A {\em bicyclic} $2k$-dimensional polytope,
$\B_{2k}(X)$,
is the convex hull of the points $\SM_{2k}(x)$, $x\in X$:
$$\B_{2k}(X)=\conv \left(\SM_{2k}(X)\right).$$
We note that  $\B_{2k}(X)$ is a 
centrally symmetric polytope
as long as one chooses $X$ to be a centrally symmetric subset of $\S^1$.
In the case of $k=2$ these polytopes were introduced and 
studied (among certain more general 
4-dimensional polytopes)
by Smilansky \cite{Smil85, Smil90}, but to the best of our knowledge
the higher-dimensional bicyclic polytopes have not yet been investigated.
Also, in \cite{Smil85} Smilansky studied the convex hull of 
$\SM_4\left(\S^1\right)$
(among convex hulls of certain more general 4-dimensional curves) but the 
convex hull of higher dimensional symmetric moment curves
has not been studied either.

We recall that a {\em face} of a convex body is the 
intersection of the body with a supporting 
hyperplane. Faces of dimension 0 are called 
{\em vertices} and faces of dimension 1 are 
called {\em edges}. Our first main result concerns 
the edges of the convex hull 
$$\B_{2k}=\conv\left(\SM_{2k}\left(\S^1\right)\right)$$
of the symmetric moment curve. Note that $\B_{2k}$ is 
centrally symmetric about the origin.

Let $\alpha \ne \beta \in \S^1$ be a pair of non-antipodal points.
 By the {\em arc with the endpoints}
$\alpha$ {\em and} $\beta$ we always mean the shorter 
of the two arcs defined by $\alpha$ and 
$\beta$.

\begin{theorem} \label{lower1} For every positive integer $k$ 
there exists a number 
$$\frac{2k-2}{2k-1}\pi\  \leq \ \psi_k  \ < \ \pi$$ 
with the following property: 
if the length of the arc with the endpoints $\alpha \ne \beta \in\S^1$ 
is less than $\psi_k$ then the interval 
$\big[\SM_{2k}(\alpha), \  \SM_{2k}(\beta)\big]$
is an edge of $\B_{2k}$ and if the length of the arc with the endpoints 
$\alpha \ne \beta \in\S^1$ is greater than $\psi_k$ 
then the interval $\big[\SM_{2k}(\alpha), \  \SM_{2k}(\beta)\big]$
is not an edge of $\B_{2k}$.
\end{theorem}
%%%%%%%%%%%%%%%%%%%%%%%%%%%%%%%%
It looks quite plausible that 
$$\psi_k=\frac{2k-2}{2k-1}$$
and, indeed, this is the case for $k=2$, cf. Section \ref{case4}.

One remarkable property of the convex hull of the trigonometric moment curve 
in $\R^{2k}$
is that it is $k$-{\em neighborly}, that is, the convex hull of any set 
of $k$ distinct points
on the curve is a $(k-1)$-dimensional face of the convex hull.
 The convex hull of the symmetric moment curve turns out to be 
{\em locally} $k$-neighborly.

\begin{theorem}\label{facegen}
For every positive integer $k$ there exists a number 
$\phi_k>0$ such that if $t_1, \ldots, t_k \in \S^1$ 
are distinct points that lie on an arc of length at most $\phi_k$, then 
$$\conv \bigl( \SM_{2k}(t_1), \ldots, \SM_{2k}(t_k) \bigr)$$
is a $(k-1)$-dimensional face of $\B_{2k}$.
\end{theorem}

From Theorems \ref{lower1} and \ref{facegen} on one hand and using
a volume trick similar to that used in \cite{LinNov} on the other hand, 
we prove 
the following results on $\fmax(d,n; j)$---{\em the maximum} 
number of $j$-faces
that a cs $d$-polytope on $n$ vertices can have.

\begin{theorem}  \label{main-thm1}
If $d$ is a fixed even number and $n\longrightarrow \infty$, then
$$
 1-\frac{1}{d-1}+o(1)  \leq
\frac{\fmax(d,n; 1)}{{n \choose 2}} \leq 1-\frac{1}{2^d} + o(1).
$$
\end{theorem}

\begin{theorem} \label{main-thm2}
If $d=2k$ is a fixed even number, $j\leq k-1$, and 
$n\longrightarrow \infty$, then
$$
c_j(d) +o(1)\leq \frac{\fmax(d,n; j)}{{n \choose j+1}} 
\leq 1-\frac{1}{2^d} + o(1),
$$
where $c_j(d)$ is a positive constant.
\end{theorem}

Some discussion is in order.
Recall that the cyclic polytope is  
{\em $\lfloor d/2 \rfloor$-neighborly}, that is, for all
$j\leq \lfloor d/2 \rfloor$, every 
$j$ vertices of $C_d(n)$ form the vertex set of a face.
Since $C_d(n)$ is a simplicial polytope, its neighborliness implies that
$f_j(C(d,n)) ={n \choose j+1}$ for $j < \lfloor d/2 \rfloor$.
Now if $P$ is a centrally symmetric polytope on $n$ vertices then no two
of its antipodal vertices are connected by an edge, and so
$f_1(P)\leq {n \choose 2}- \frac{n}{2}$.
In fact, as was recently shown by 
Linial and the second author \cite{LinNov},
this inequality is strict as long as $n>2^d$. This leads one to wonder
how big the gap between $\fmax(d,n;1)$ and ${n \choose 2}$ is and whether
$\fmax(d,n; j)$, for $j < \lfloor d/2 \rfloor$,
is on the order of $n^{j+1}$.
Theorems \ref{main-thm1} and \ref{main-thm2} (see also
Propositions \ref{UB1} and \ref{UB2} below)
provide (partial) answers to those questions. 

Let us fix an even dimension $d=2k$ and let 
$X \subset \S^1$ be a set of $n$ equally spaced 
points, where $n$ is an even number. 
We conjecture that for every integer $j \leq k-1$ 
$$\limsup_{n \longrightarrow +\infty} 
\frac{f_j\left(\B_{2k}(X)\right)}{{n \choose j+1}}=
\limsup_{n \longrightarrow +\infty} \frac{\fmax(d,n;j)}{{n \choose j+1}}.$$

It is also worth mentioning that recently there has been a lot of interest
in the problems surrounding neighborliness and face numbers 
of cs polytopes in connection
to statistics and error-correcting codes, 
see \cite{Dono1, Dono2, DonoTan, RudelVersh}.
In particular, it was proved in \cite{DonoTan} that 
for large $n$ and $d$, if $j$ is bigger than a certain threshold value,
then the ratio between 
the expected number of $j$-faces of 
a random cs $d$-polytope with $n$ vertices
and ${n \choose j+1}$ is
smaller than $1-\epsilon$ for some positive constant $\epsilon$.
The upper bound part of Theorem~\ref{main-thm2} 
provides a real reason for this phenomenon:
the expected number of $j$-faces is ``small'' 
because the $j$th face number of {\em every} cs polytope is ``small".

The structure of the paper is as follows.
In Section 2 we prove the upper bound parts of 
Theorems \ref{main-thm1} and \ref{main-thm2}.
In Section 3 we discuss bicyclic polytopes and their relationship 
to non-negative trigonometric polynomials and self-inversive polynomials.
Section 4 contains new short proofs of results
originally due to Smilansky
on the faces of 4-dimensional bicyclic polytopes.
It serves as a warm-up for Sections 5 and 6  in which we
prove Theorems \ref{lower1} and \ref{facegen} as well as the 
lower bound parts of 
Theorems \ref{main-thm1} and \ref{main-thm2}. 
We discuss 2-faces of $\B_6$ as well as $\fmax(2k, n; j)$ for $j \geq k$
in Section 7, where we also state several open questions.

\section{Upper bounds on the face numbers}
The goal of this section is to prove the upper bound parts of 
Theorems \ref{main-thm1} and \ref{main-thm2}. The proof uses
a volume trick similar to the one utilized in the proof of
the Danzer-Gr\"unbaum theorem on 
the number of vertices of antipodal polytopes \cite{DanzGr} 
and more recently in \cite[Theorem 1]{LinNov}, where it was used
to estimate maximal possible neighborliness of cs polytopes.

The upper bound part of Theorem \ref{main-thm1}
is an immediate consequence of the following more precise result.
\begin{proposition}  \label{UB1}
Let $P\subset \R^d$ be a cs $d$-polytope on $n$ vertices.
Then 
$$f_1(P) \leq \frac{n^2}{2} \left(1-2^{-d} \right).$$
\end{proposition}
\proof
Let $V$ be the set of vertices of $P$.
For every vertex $u$ of $P$ we define
$$P_u:=P+u \subset 2P$$ to be a translate of $P$, 
where ``+" denotes the Minkowski addition. 
We claim that if the polytopes $P_u$ and $P_v$ 
have intersecting interiors then the vertices $u$ and $-v$
are not connected by an edge. (Note that this includes the case of $u=v$,
since clearly $\Int(P_v)\cap\Int(P_v)\neq\emptyset$ and $(v, -v)$
is not an edge of $P$.)  Indeed, the assumption
$\Int(P_u)\cap\Int(P_v)\neq \emptyset$ implies that there exist 
$x,y\in\Int(P)$ such that $x+u=y+v$, or equivalently, 
that $(y-x)/2=(u-v)/2$. Since $P$ is centrally symmetric, and 
$x,y\in\Int(P)$, the point $q:=(y-x)/2$ is an interior point of $P$.
As $q$ is also the barycenter of the line segment connecting
$u$ and $-v$, this line segment is not an edge of $P$.

Let us normalize the Lebesgue measure $dx$ in $\R^d$ 
in such a way that $\vol(2P)=1$ and hence
$$\vol(P)=\vol\left(P_u\right)=2^{-d} \quad \text{for all} \quad u \in V.$$
For a set $A \subset \R^d$, let $[A]: \R^d \longrightarrow \R$ be the 
indicator of $A$, that is,
$[A](x) =1$ for $x \in A$ and $[A](x)=0$ and let us define
$$h=\sum_{u \in V} [\Int P_{u}].$$
Then
$$\int_{2P} h \ dx =n2^{-d},$$
and hence by the H\"older inequality
$$\int_{2P} h^2 \ dx  \geq n^2 2^{-2d}.$$
On the other hand, the first paragraph of the proof implies that
$$\int_{2P} h^2(x) \ dx 
=\sum_{u,v \in V} \vol\left(P_u \cap P_v \right) \leq n2^{-d} +2^{-d+1}
\left({n \choose 2}-f_1(P) \right),$$ 
and the statement follows.
\endproof 

As a corollary, we obtain the following upper bound on $f_j(P)$  
for any $1 \leq j \leq d/2$. This upper bound implies the upper bound part
of Theorem \ref{main-thm2}.

\begin{proposition}    \label{UB2} 
Let $P\subset \R^d$ be a cs $d$-polytope with $n$ vertices,
and let $j\leq (d-2)/2$.
Then 
$$ f_j(P) \leq \frac{n}{n-1}\left(1-2^{-d}\right) {n \choose j+1}.$$
\end{proposition}

\proof We rely on Proposition \ref{UB1} and two additional results. 

The first result is an adaptation
of the well-known perturbation argument, 
see, for example, \cite[Section 5.2]{Grunb},
to the centrally symmetric situation. 
Namely, we claim that for every cs $d$-polytope $P$ 
there exists a simplicial cs $d$-polytope $Q$ such that
$f_0(P)=f_0(Q)$  and $f_j(P)\leq f_j(Q)$ for all $1\leq j\leq d-1$. 
The polytope $Q$ is obtained 
from $P$ by pulling the vertices of $P$ in a generic way but so as 
to preserve the symmetry.
The proof is completely similar to that of  \cite[Section 5.2]{Grunb} 
and hence is omitted.

The second result states that for any $(d-1)$-dimensional simplicial 
complex $K$ with $n$ vertices, we have 
$$f_j(K)\leq f_1(K) {n \choose j+1} / {n \choose 2} 
\quad \text{for} \quad 1 \leq j \leq d-1.$$
The standard double-counting argument goes as follows: 
every $j$-dimensional simplex of $K$
contains exactly ${j+1 \choose 2}$ edges and every edge of $K$ is contained in 
at most ${n-2 \choose j-1}$ of the $j$-dimensional simplices of $K$. Hence
$$f_j(K)/f_1(K) \leq {n-2 \choose j-1}/{j+1 \choose 2}
={n \choose j+1}/{n \choose 2}.$$

The statement now follows by Proposition \ref{UB1}.
\endproof

\section{Faces and polynomials}

In this section, we relate the facial structure of the 
convex hull $\B_{2k}$ of the symmetric 
moment curve (Section \ref{curve}) to properties of trigonometric 
and complex polynomials from particular families.

\subsection{Preliminaries}\label{prelim}

A proper face of a convex body $B\subset\R^{2k}$ 
is the intersection of $B$ with its supporting hyperplane, 
that is,  the intersection of $B$
with the zero-set of an affine function
$$A(x)=\alpha_0+ \alpha_1\xi_1+ \ldots+ \alpha_{2k}\xi_{2k} \quad \text{for}
\quad x=\left(\xi_1, \ldots, \xi_{2k}\right)$$ which satisfies $A(x)\geq 0$ for
all $x \in B$.

A useful observation is that $\B_{2k}$ remains invariant under 
a one-parametric group
of rotations that acts transitively on $\SM_{2k}\left(\S^1\right)$. 
Such a rotation is 
represented by the $2k \times 2k$ block-diagonal matrix with the 
$j$th block being
$$\left(\begin{matrix} 
\cos (2j-1)\tau & \sin(2j-1)\tau \\
 -\sin(2j-1)\tau & \cos (2j-1) \tau \end{matrix}
\right)$$
for $\tau \in \R$. 
If $\{t_1, \ldots, t_s\} \subset \S^1$ are distinct points such that
$$\conv\Bigl( \SM_{2k}(t_1), \ldots, \SM_{2k}(t_s) \Bigr)$$
is a face of $\B_{2k}$ and the points $t_1', \ldots, t_s' \in \S^1$ 
are obtained from
$t_1, \ldots, t_s$ by a rotation $t_i'=t_i +\tau$ for $i=1, \ldots, s$ 
of  $\S^1$, then 
$$\conv\Bigl( \SM_{2k}(t_1'), \ldots, \SM_{2k}(t_s') \Bigr)$$
is a face of $\B_{2k}$ as well.

Finally, we note that the natural projection 
$\R^{2k} \longrightarrow \R^{2k'}$ for $k'<k$ that 
erases the last $2k-2k'$ coordinates maps $\B_{2k}$ 
onto $\B_{2k'}$ and $\B_{2k}(X)$
onto $\B_{2k'}(X)$  Hence, if for some 
sets $Y \subset X \subset \S^1$ the set
$\conv\left(\SM_{2k'}(Y)\right)$ is a face of 
$\B_{2k'}(X)$, then  $\conv\left(\SM_{2k}(Y)\right)$ 
is a face of $\B_{2k}(X)$.
 
\subsection{Raked trigonometric polynomials}
The value of an affine function $A(x)$ on the symmetric moment curve 
$\SM_{2k}$ is represented by a trigonometric polynomial 
\begin{equation}  \label{raked-trig-poly}
A(t)=c+ \sum_{j=1}^k a_j \cos (2j-1)t+ 
              \sum_{j=1}^k b_j\sin (2j-1)t .\end{equation}
Note that all summands involving the even terms $\sin 2jt$ and $\cos 2jt$
except for the constant term vanish from $A(t)$. 
We refer to such trigonometric polynomials as 
{\em raked trigonometric polynomials} of degree $2k-1$. 
As before, it is convenient to think 
of $A(t)$ as defined on $\S^1=\R / 2\pi \Z$. 

Admitting, for convenience, the whole body $\B_{2k}$ and the empty set 
as faces of $\B_{2k}$, we obtain
the following result.

\begin{lemma}  \label{faces-descr1}
The faces of $\B_{2k}$ are defined
by the raked trigonometric polynomials of degree $2k-1$ 
that are non-negative on $\S^1$.
If $A(t)$ is such a polynomial and 
$\{t_1,\ldots, t_s\} \subset \S^1$ is the set of its zeroes, 
then the face of $\B_{2k}$ defined by $A(t)$ is the convex hull of 
$\bigl\{\SM_{2k}(t_1), \ldots, \SM_{2k}(t_s)\bigr\}$.
\end{lemma}

\subsection{Raked self-inversive polynomials}\label{raked-rel}

Let us substitute $z=e^{it}$ in equation (\ref{raked-trig-poly}). 
Using that
$$\cos (2j-1)t=\frac{z^{2j-1} +z^{1-2j}}{2} 
\quad \text{and} \quad \sin(2j-1)t =\frac{z^{2j-1}-z^{1-2j}}{2i}$$
we can write $A(t)=z^{-2k+1} D(z)$, where 
$$D(z)=c z^{2k-1}+\sum_{j=1}^k \frac{a_j -i b_j}{2} z^{2j+2k-2} 
+ \sum_{j=1}^k \frac{a_j+i b_j}{2} z^{2k-2j}.$$
In other words, $D(z)$ is a polynomial satisfying  
\begin{equation}\label{selfinv}
D(z)=z^m \overline{D\bigl(1/\overline{z}\bigr)},  \quad \mbox{where }
m=4k-2
\end{equation}
and such that 
\begin{equation}\label{raked}
D(z)=cz^{2k-1} + \sum_{j=0}^{2k-1} d_{2j} z^{2j},
\end{equation}
so that all odd terms with the possible exception of the middle term vanish.

The polynomials $D(z)$ satisfying equation (\ref{selfinv}) are well studied 
and  known in the literature by the name  
{\em self-inversive polynomials} (see 
for instance \cite[Chapter 7]{Small}).
In analogy with raked trigonometric polynomials, 
we refer to polynomials $D$ satisfying both
equations (\ref{selfinv}) and (\ref{raked}) 
as {\em raked self-inversive polynomials} of degree $m$.
We note that any polynomial $D(z)$ satisfying 
(\ref{selfinv}) with $m=4k-2$ and 
(\ref{raked}) gives rise to a raked trigonometric polynomial $A(t)$ such that 
$A(t)=z^{-2k+1}D(z)$ for $z=e^{it}$.
Hence we obtain the following  restatement of Lemma \ref{faces-descr1}.

\begin{lemma}   \label{faces-descr2}
The faces of $\B_{2k}$
are defined by the raked self-inversive polynomials of degree $4k-2$
all of whose roots of modulus one have even multiplicities.
If $D(z)$ is such a polynomial and $\bigl\{e^{it_1}, \ldots, e^{it_s}\bigr\}$
is the set of its roots of modulus 1, then the face of $\B_{2k}$
defined by $D(z)$ is the convex hull of 
$\bigl\{\SM_{2k}(t_1), \ldots, \SM_{2k}(t_s)\bigr\}$.
\end{lemma}

Let $D(z)$ be a polynomial satisfying equation (\ref{selfinv}), 
and let 
$$M=\bigl\{ \zeta_1, \ldots, \zeta_1, \zeta_2, \ldots, \zeta_2, \ldots, 
\zeta_s, \ldots, \zeta_s \bigr\}$$
be the multiset of all roots of $D$ where each root is listed the
 number of times equal to its 
multiplicity. We note that if $\deg D=m$ then $0 \notin M$ and $|M|=m$. 
We need a straightforward characterization of the raked self-inversive
polynomials in terms of their zero multisets $M$.

\begin{lemma}  \label{roots-selfinv}
A multiset  $M \subset \C$ of size $|M|=4k-2$ is the multiset of 
roots of a raked self-inversive polynomial 
of degree $4k-2$ if and only if  
$$\overline{M}=M^{-1},$$
that is, $\zeta \in M$ if and only if $\overline{\zeta}^{-1} \in M$ 
and the multiplicities of 
$\zeta$ and $\overline{\zeta}^{-1}$ in $M$ are equal, 
and 
$$\sum_{\zeta \in M} \zeta^{2j-1}=0 \quad \text{for} \quad j=1, \ldots, k-1.$$
\end{lemma}
\proof It is known and not hard to see
that $M$ is the zero-multiset of 
a self-inversive polynomial if and only if $M^{-1}=\overline{M}$
\cite[p.~149, p.~228]{Small}. Indeed, if $D(0) \ne 0$ then (\ref{selfinv}) 
implies 
$\overline{M}=M^{-1}$. Conversely, suppose that  
$M$ is the multiset satisfying $\overline{M}=M^{-1}$.
Then 
$$\prod_{\zeta \in M} |\zeta|=1,$$
and hence we can choose numbers $a_{\zeta}$ such that 
$$\prod_{\zeta \in M} {\overline a_{\zeta} \over a_{\zeta}}= 
\prod_{\zeta \in M} (-\zeta).$$
Then the polynomial
$$D(z)=  \prod_{\zeta \in M} a_{\zeta} (z -\zeta)$$
satisfies (\ref{selfinv}) with $m=|M|$.

Let 
$$s_p=\sum_{\zeta \in M} \zeta^p \quad \text{and let} 
\quad D(z)=\sum_{p=0}^m d_p z^p.$$
Using Newton's formulas to express elementary symmetric functions 
in terms of power sums, we get
$$pd_{m-p}+\sum_{j=1}^p s_j d_{m-p+j}=0 \quad \text{for} \quad p=1,2, 
\ldots, m.$$
Since $m=4k-2$ is even, we conclude that 
$$d_1=d_3 =\ldots = d_{2k-3}=0 \quad \text{if and only if} 
\quad s_1=s_3 = \ldots= s_{2k-3}=0,$$
which completes the proof.
\endproof

We conclude this section with the description of a
particular family of faces of $\B_{2k}$.

\subsection{Simplicial faces of $\B_{2k}$} \label{simplicial}

Let 
$$A(t)=1-\cos\bigl((2k-1)t \bigr).$$
Clearly, $A(t)$ is a raked trigonometric polynomial and 
$A(t) \geq 0$ for all $t \in \S^1$.
Moreover, $A(t)=0$ at the $2k-1$ points 
$$\tau_j={2 \pi j \over 2k-1} \quad \text{for} \quad j=1, \ldots, 2k-1$$
on the circle $\S^1$, which form the vertex set of a regular $(2k-1)$-gon.
By Lemma \ref{faces-descr1} the set 
$$\Delta_0=\conv \Bigl(\SM_{2k}(\tau_1), \ldots, \SM_{2k}(\tau_{2k-1}) \Bigr)$$
is a face of $\B_{2k}$. 

One can observe that $\Delta_0$ is a $(2k-2)$-dimensional regular simplex, 
since the cyclic permutation of the vertices 
$$\tau_1 \longmapsto \tau_2 \longmapsto \ldots \longmapsto \tau_{2k-1} 
\longmapsto \tau_1$$
gives rise to an orthogonal transformation of $\R^{2k}$ which maps 
$\Delta_0$ onto itself
and also maps $\B_{2k}$ onto itself, cf. Section \ref{prelim}.
Furthermore, we have a one-parametric family of 
simplicial faces 
$$\Delta_{\phi} =\conv\Bigl(\SM_{2k}\left(\tau_1 + \tau \right), 
\ldots, \SM_{2k}\left(\tau_{2k-1} +\tau 
\right) \Bigr) \quad \text{for} \quad 0 \leq \tau < 2 \pi$$
of $\B_{2k}$. The dimension of the boundary of $\B_{2k}$ is $(2k-1)$,
so this one-parametric family of simplices covers a ``chunk'' 
of the boundary of $\B_{2k}$. 

\section{The faces of $\B_4$}\label{case4}

In this section we provide a complete characterization 
of the faces of $\B_4$. 
This result is not new, it was proved by Smilansky \cite{Smil85} who 
also described the facial structure of the convex hull
of the more general curve $(\cos pt, \sin pt, \cos qt, \sin qt)$, where
$p$ and $q$ are any positive integers. Our proof serves as a warm-up
for the following section where we discuss the edges of $\B_{2k}$ for $k>2$.

\begin{theorem}\cite{Smil85}    \label{Smil-thm}
The proper faces of $\B_4$ are

\begin{itemize}
\item[(0)] The 0-dimensional faces (vertices)
$$\SM_4(t), \quad t \in \S^1;$$

\item[(1)]  The 1-dimensional faces (edges)
$$[\SM_4(t_1),\ \SM_4(t_2)],$$ where 
$t_1 \ne t_2$ are the endpoints of an arc of $\S^1$ of length less 
than $2\pi/3$; and

\item[(2)] The 2-dimensional faces (equilateral triangles)
$$\Delta_t=\conv\Bigl(\SM_4(t),\ \SM_4(t+2\pi/3),\ \SM_4(t+4\pi/3)\Bigr), 
\quad t \in \S^1.$$ 

\end{itemize}
\end{theorem}

\proof
We use Lemma \ref{faces-descr2}. 
A face of $\B_4$ is determined by a raked self-inversive 
polynomial $D$ of degree 6. Such a polynomial $D$ 
has at most 3 roots on the circle 
$\S^1$, each having an even multiplicity. 
Furthermore, by Lemma \ref{roots-selfinv}, 
the sum of all the roots of $D$ is 0.
Therefore, we have the following three cases.

Polynomial $D$ has 3 double roots $\zeta_1=e^{it_1}, 
\zeta_2=e^{it_2}, \zeta_3=e^{it_3}$. 
Since $\zeta_1+\zeta_2+\zeta_3=0$, the points $t_1, t_2,t_3 \in \S^1$ 
form the vertex set of 
an equilateral triangle, and we obtain the 2-dimensional 
face defined in Part (2).

Polynomial $D$ has two double roots $\zeta_1=e^{it_1}, 
\zeta_2=e^{i t_2}$, and a pair 
of simple roots $\zeta$ and $\overline{\zeta}^{-1}$ with $|\zeta| \ne 1$. 
Applying a rotation,
if necessary, we may assume without loss of generality that $t_1=-t_2=t$. 
Since we must have 
$$\zeta + \overline{\zeta}^{-1} + 2 e^{it} + 2 e^{-it} =0,$$
we conclude that $\zeta \in \R$.  Hence the equation reads
$$\zeta+ \zeta^{-1} =- 4 \cos t \quad \text{for some} 
\quad \zeta \in \R, \quad |\zeta| \ne 1.$$
If $|\cos t| >1/2$ then the solutions $\zeta, \zeta^{-1}$ 
of this equation are indeed real and satisfy 
 $|\zeta|, |\zeta^{-1}| \ne 1$. If $|\cos t| \leq 1/2$
then the solutions $\zeta, \zeta^{-1}$ 
form a pair of complex conjugate numbers satisfying
$|\zeta|=|\zeta^{-1}|=1$. Therefore, the interval 
$[\SM_{4}(-t), \ \SM_4(t)]$ is a face of $\B_4$ if and only 
if $-\pi/3<t< \pi/3$ or $2\pi/3 <t < 4\pi/3$, 
so we obtain the 1-dimensional faces as in Part (1).

Finally, we conclude that since $\B_4$ 
must have at least one 0-dimensional face
(vertex)  and that the vertices of $\B_4$ may only be of the type $\SM_4(t)$
for some $t \in \S^1$, one of the points $\SM_4(t)$ must be a vertex of $\B_4$.
Because of rotational invariance (cf. Section \ref{prelim}), 
all the points $\SM_4(t)$ 
are the vertices
of $\B_4$, which concludes the proof.
\endproof

\section{Edges of $\B_{2k}$}

In this section we prove Theorems \ref{lower1} and  
\ref{main-thm1}. Our main tool is a certain deformation of simplicial 
faces of $\B_{2k}$, cf. Section \ref{simplicial}.

\subsection{Deformation}\label{deformation}

Let $M$ be a finite multiset of non-zero complex numbers 
such that $M=M^{-1}$. In other words,
for every $\zeta \in M$ we have $\zeta^{-1} \in M$ 
and the multiplicities of $\zeta$ and $\zeta^{-1}$ 
in $M$ are equal. In addition, we assume that the multiplicities of 1
 and $-1$ in $M$ are even, 
possibly 0. For every $\lambda \in \R \setminus \{ 0\}$ we define the multiset 
$M_{\lambda}$, which we call a {\em deformation} of $M$, as follows.

We think of $M$ as a multiset of unordered pairs $\{\zeta, \zeta^{-1}\}$. 
For every such pair, we
consider the equation

\begin{equation}\label{deform}
z+z^{-1} = \lambda \left(\zeta + \zeta^{-1}\right).
\end{equation}

We let $M_{\lambda}$ to be the multiset consisting of the pairs 
$\{z, z^{-1} \}$ of solutions of (\ref{deform}) 
as $\{\zeta, \zeta^{-1}\}$ range over $M$. Clearly, $|M_{\lambda}|=|M|$ and 
$M_{\lambda}^{-1}=M_{\lambda}$. In addition, if $\overline{M}=M$ then 
$\overline{M_{\lambda}}=M_{\lambda}$, since $\lambda$ in (\ref{deform}) 
is real.

%\bigskip

Our interest in the deformation $M \longmapsto M_{\lambda}$ 
is explained by the following lemma.

\begin{lemma}\label{deformpoly} Let $D(z)$ be a raked self-inversive 
polynomial of degree $4k-2$ 
with real coefficients and such that $D(0) \ne 0$. 
Let $M$ be the multiset of the roots of $D$ and 
suppose that both $1$ and $-1$ have an even, 
possibly 0, multiplicity in $M$.  
Then, for any real $\lambda \ne 0$, 
the defomation $M_{\lambda}$ of $M$
is the multiset of the roots of a raked self-inversive polynomial 
$D_{\lambda}(z)$ of degree 
$4k-2$ with real coefficients.
\end{lemma}
\proof We use Lemma \ref{roots-selfinv}.
 Since $D$ has real coefficients, we have 
$M=\overline{M}$, so by Lemma~\ref{roots-selfinv}, 
we have $M=M^{-1}$ as well.
Clearly, $M_{\lambda}=M_{\lambda}^{-1}$ and 
${M_{\lambda}}=\overline{M_{\lambda}}$,
so $M_{\lambda}$ is the multiset of the roots 
of a self-inversive real polynomial $D_{\lambda}$ 
of degree $4k-2$. It remains to check that
$$\sum_{\zeta \in M_{\lambda}} \zeta^{2j-1} =0 
\quad \text{for} \quad j=1, \ldots, k-1.$$
We have 
\begin{equation}\label{binom}
\left(x+x^{-1}\right)^{2n-1}=\sum_{m=1}^n {2n-1 \choose n+m-1} 
\left(x^{2m-1} + x^{-2m+1} \right).
\end{equation}
Since by Lemma \ref{roots-selfinv}
$$\sum_{\zeta \in M} \zeta^{2j-1}=\sum_{\zeta \in M} \zeta^{1-2j} =0 
\quad \ \text{for} \quad j=1, \ldots, k-1,$$
it follows by formula (\ref{binom}) that
$$\sum_{\zeta \in M} \left(\zeta + \zeta^{-1}\right)^{2j-1}=0 
\quad \text{for} \quad j=1, \ldots, k-1.$$
Therefore, by (\ref{deform}), we have
$$\sum_{\zeta \in M_{\lambda}} 
\left( \zeta + \zeta^{-1}\right)^{2j-1}=0 \quad \text{for} \quad 
j=1, \ldots, k-1,$$
from which by (\ref{binom}) we obtain 
$$\sum_{\zeta \in M_{\lambda}} \zeta^{2j-1}=
\frac{1}{2}\sum_{\zeta \in M_{\lambda}} \left(\zeta^{2j-1} 
+\zeta^{-2j+1} \right)=0
\quad \text{for} \quad j=1, \ldots, k-1,$$
as claimed. Hence $D_{\lambda}(z)$ is a raked self-inversive polynomial.
\endproof

To prove Theorem \ref{lower1} we need another auxiliary result.
%We use Lemma \ref{deformpoly} to prove Theorem \ref{lower1}.

\begin{lemma}\label{long-short}
Let $\alpha, \beta \in \S^1$ be such that the interval 
$[\SM_{2k}(\alpha),\ \SM_{2k}(\beta)]$ is an edge of $\B_{2k}$ 
and let $\alpha' \ne \beta' \in \S^1$ 
be some other points such that the arc with the endpoints 
$\alpha', \beta'$ is shorter than the 
arc with the endpoints $\alpha, \beta \in \S^1$. Then the interval 
$[\SM_{2k}(\alpha'),\ \SM_{2k}(\beta')]$ is an edge of $\B_{2k}$.
\end{lemma}
\proof Because of rotational invariance, we assume,  without loss 
of generality, that 
$\alpha=\tau$ and $\beta=-\tau$ for some 
$0< \tau < \pi/2$. Let $A(t)$ be a raked trigonometric
polynomial that defines the edge $[\SM_{2k}(\alpha),\ \SM_{2k}(\beta)]$, 
see Lemma  \ref{faces-descr1}. 
Hence $A(t )\geq 0$ for all $t \in \S^1$ and $A(t)=0$ if and only 
$t =\pm \tau$. Let $A_1(t)=A(t)+A(-t)$. 
Then $A_1(t)$ is a raked trigonometric polynomial 
such that $A_1(t) \geq 0$ for all $t \in \S^1$ and $A_1(t)=0$ 
if and only if $t =\pm \tau$.
Furthermore, we can write
$$A_1(t)=c+\sum_{j=1}^k a_j \cos (2j-1)t $$
for some real $a_j$ and $c$. 
Moreover,  we assume, without loss of generality, that $a_k \ne 0$. 
(Otherwise choose $k'$ to be the largest index $j$ with $a_j  \ne 0$ and 
 project $\B_{2k}$ onto $\B_{2k'}$, cf. Section \ref{prelim}.)
Hence the polynomial $D(z)$
defined by $A_1(t)=z^{-2k+1}D(z)$ for $z=e^{it}$, 
see Section \ref{raked-rel}, is a raked self-inversive
polynomial of degree $4k-2$ with real coefficients 
satisfying $D(0) \ne 0$. Moreover,
the only roots of $D(z)$ that lie on the circle $|z|=1$ are 
$e^{i\tau}$ and $e^{-i\tau}$ and those roots 
have equal even multiplicities.

Let us choose an arbitrary $0< \tau' <\tau$ and let
$$\lambda=\frac{\cos \tau'}{\cos \tau}>1.$$
Let $D_{\lambda}$ be the raked self-inversive polynomial of degree $4k-2$ 
whose existence is established by Lemma \ref{deformpoly}. Since
$$e^{i \tau'}+e^{-i \tau'}=\lambda\left(e^{i \tau} +e^{-i \tau}\right),$$
the numbers $e^{i \tau'}$ and $e^{-i \tau'}$ are roots of 
$D_{\lambda}$ of even multiplicity.
Moreover, suppose that $z$ is a root of  $D_{\lambda}$ such that  $|z|=1$. 
Then $z+z^{-1} \in \R$ and $-2 \leq z+z^{-1} \leq 2$. By 
(\ref{deform}) it follows that there is a pair $\zeta, \zeta^{-1}$ 
of roots of $D$ such that
$\zeta +\zeta^{-1} \in \R$ and $-2 <|\zeta +\zeta^{-1}| <2$. 
It follows then that $|\zeta|=|\zeta^{-1}|=1$,
from which, necessarily, 
$\left\{\zeta, \zeta^{-1}\right\}=\left\{e^{i\tau},\ e^{-i\tau} \right\}$
and hence $\left\{z, z^{-1}\right\}=\left\{e^{i \tau'}, e^{-i \tau'}\right\}$.
Therefore, by Lemma 
\ref{faces-descr2}, $[\SM_{2k}(-\tau'), \ \SM_{2k}(\tau') ]$ is an edge of 
$\B_{2k}$. 
By rotational invariance, it follows that 
$[\SM_{2k}(\alpha'), \ \SM_{2k}(\beta')]$ is 
an edge of $\B_{2k}$, where points 
$\alpha', \beta'$ are obtained from $\tau', -\tau'$ by 
a rotation of $\S^1$.
\endproof

We are now ready to complete the proof of Theorem \ref{lower1}.

\smallskip\noindent {\it Proof of Theorem \ref{lower1}: \ } 
In view of Lemma \ref{long-short},
it remains to show that one can find an arbitrarily small 
$\delta>0$ and two points 
$\alpha, \beta \in \S^1$ such that
the interval $[\SM_{2k}(\alpha), \SM_{2k}(\beta)]$ 
is an edge of $\B_{2k}$ and 
the length of the arc with the endpoints $\alpha$ and $\beta$ is at least 
$\displaystyle \frac{2 \pi (k-1)}{2k-1}-\delta$.

 Let us consider the polynomial
$$D(z)=\left(z^{2k-1}-1\right)^2=z^{4k-2} -2 z^{2k-1}+1.$$
Clearly, $D(z)$ is a raked self-inversive polynomial of degree 
$4k-2$ and the multiset $M$ of the 
roots of $D$ consists of all roots of unity of degree $2k-1$, 
each with multiplicity 2.
In fact, $D(z)$ defines a simplicial face of $\B_{2k}$, 
cf. Section \ref{simplicial}. For 
$\epsilon>0$ let us consider the deformation $D_{1+\epsilon}(z)$ 
of $D(z)$ and its roots,
see Lemma \ref{deformpoly}. 

In view of equation (\ref{deform}), for all sufficiently 
small $\epsilon>0$, the multiset $M_{1+\epsilon}$ 
of the roots of $D_{1+\epsilon}$ consists of two 
positive simple real roots found from the equation
$$z+z^{-1}=2(1+\epsilon),$$
that are the deformations of the double root at 1 
(one of them is greater than and the other
is less than 1), and $2k-2$ double roots on the unit circle, 
found from the equation
$$z+z^{-1}=2(1+\epsilon) \cos \frac{2 \pi j}{2k-1} \quad \text{for} 
\quad j=1, \ldots, 2k-2,$$
that are the deformations of the remaining $2k-2$ roots of unity.

Let $\zeta_j(\epsilon)$ be the deformation of the root 
$$\zeta_j =\cos \frac{2 \pi j}{2k-1} +i \sin \frac{2 \pi j}{2k-1} 
\quad \text{for}\quad j=1, \ldots, 2k-2$$
that lies close to $\zeta_j$ if $\epsilon>0$ is small enough, see Figure 1.
\begin{figure}
$$\epsffile{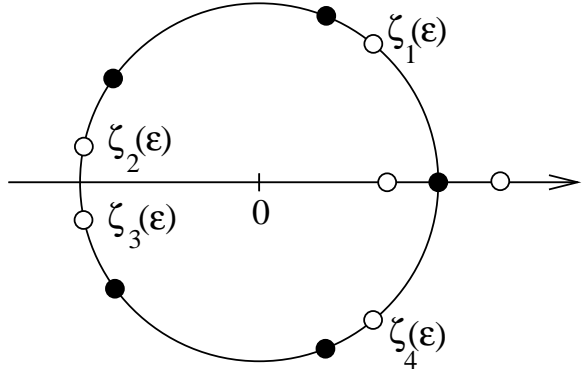}$$
\caption{The roots of unity (black dots) and their deformations
(white dots) for $k=3$}
\end{figure}
Thus we have 
\begin{equation}\label{conjugate}
\zeta_j^{-1}(\epsilon)=\overline{\zeta_j}(\epsilon)=\zeta_{2k-1-j}(\epsilon).
\end{equation}
Let 
$$\zeta_j(\epsilon)=e^{i \alpha_j}
\quad \text{where} \quad 0 < \alpha_j < 2 \pi \quad
 \text{for} \quad j=1, \ldots, 2k-2.$$
Then 
$$\cos \alpha_j = (1+\epsilon)\cos \frac{2 \pi j}{2k-1},$$
and hence 
\begin{equation}\label{infinitesimal}
\alpha_j= \frac{2 \pi j}{2k-1}-
 \epsilon \ \ctg  \frac{2 \pi j}{2k-1} +O(\epsilon^2).
\end{equation}

Let us prove that the interval 
$$\left[ \SM_{2k}(\alpha_1), \ \SM_{2k}(\alpha_k) \right]$$
is an edge of $\B_{2k}$. 
We obtain this edge as the intersection of two faces of $\B_{2k}$.

By Lemma~\ref{faces-descr2}, 
\begin{equation}\label{firstface}
\conv \Bigl( \SM_{2k}(\alpha_1), \ldots, \SM_{2k}(\alpha_{2k-2}) \Bigr)
\end{equation}
is a face of $\B_{2k}$. 
 The second face is obtained by a rotation of (\ref{firstface}).
Namely, let us consider the clockwise rotation of the circle $|z|=1$ 
which maps 
$\zeta_{k-1}(\epsilon)$ onto $\zeta_1(\epsilon)$. Because of (\ref{conjugate}) 
this rotation also maps $\zeta_{2k-2}(\epsilon)$ onto 
$\zeta_k(\epsilon)$. Furthermore, for $j=1, \ldots, 2k-2$ let us define 
 $$\zeta_j'(\epsilon)=e^{i \alpha_j'} \quad \text{for} \quad 0 < \alpha_j' 
< 2\pi$$
as the image of $\zeta_{j+k-2}(\epsilon)$ if  $j \leq k$, as the image of 
$\zeta_{j-k-1}(\epsilon)$ 
if $j >k+1$, and as the image of $\zeta_{k-2}(\epsilon)$ if $j=k+1$
under this rotation. Using  (\ref{infinitesimal}), we conclude that 

\begin{eqnarray*}
 \alpha_j'-\alpha_j & = &
\epsilon\left(\ctg \frac{2 \pi j}{2k-1} -
\ctg \frac{2\pi j -3\pi}{2k-1} -\ctg \frac{\pi}{2k-1} -\ctg \frac{2\pi}{2k-1}
\right) \\ 
&+ & O\left(\epsilon^2 \right) \quad \text{for} 
\quad 1 \leq j \leq 2k-2, \quad j \ne k+1 \end{eqnarray*}
and
$$\alpha_{k+1}'=o(1) \quad \text{as} \quad \epsilon \longrightarrow 0+.$$
If $\epsilon>0$ is sufficiently small then  
for $j \ne k+1$ the value of $\alpha_j'$ is close to $\alpha_j$ and strictly 
smaller than $\alpha_j$ unless $j=1$ or $j=k$, in which case
 the two values are equal. 
Furthermore,
$\alpha'_{k+1} \ne \alpha_j$ for any $j$.
By rotational invariance, see Section \ref{prelim},
\begin{equation}\label{secondface}
\conv \Bigl( \SM_{2k}(\alpha_1'), \ldots, \SM_{2k}(\alpha_{2k-2}') \Bigr)
\end{equation}
is a face of $\B_{2k}$ as well. 
Since faces (\ref{firstface}) and (\ref{secondface}) intersect 
along the interval 
$$\left[\SM_{2k}(\alpha_1) \ \SM_{2k}(\alpha_k)\right],$$ 
this interval is  an edge 
of $\B_{2k}$. Since $\alpha_1$ and $\alpha_k$ are the endpoints of 
an arc of length 
$$\pi \frac{2k-2}{2k-1} -O(\epsilon),$$
the statement follows.  
\endproof

\smallskip\noindent {\it Proof of Theorem \ref{main-thm1}: \ } 
The upper bound follows
by Proposition \ref{UB1}. To prove the lower bound, 
let us consider the polytope
$\B_{2k}(X)$, where $X \subset \S^1$ is the set of $n$ 
equally spaced points ($n$ is even).
The lower bound follows by Theorem \ref{lower1}.
\endproof

\section{Faces of $\B_{2k}$}

In this section, we prove Theorems \ref{facegen} and \ref{main-thm2}. 
Theorem \ref{facegen} is deduced from the following proposition.

\begin{proposition} \label{rakegen}
For every positive integer $k$ there exists a number $\phi_k>0$ 
such that for every set of 
$2k$ distinct points $t_1, \ldots, t_{2k} \in \S^1$ lying on 
an arc of length at most $\phi_k$
there exists a raked trigonometric polynomial 
$$A(t)=c_0 + \sum_{j=1}^k a_j \sin (2j-1)t +\sum_{j=1}^k b_j \cos (2j-1)t$$
such that $A(t)=0$ for $t \in \S^1$ if and only if $t=t_j$ for some $j=1,
 \ldots, 2k$.
\end{proposition}

To prove Proposition \ref{rakegen},  we establish first that the curve
$\SM_{2k}(t)$ is nowhere locally flat.

\begin{lemma} \label{noflat}
Let 
$$\SM_{2k}(t)=\Bigl(\cos t, \ \sin t, \ \cos 3t, \ \sin 3t, \ \ldots,\  
\cos (2k-1)t,\ \sin (2k-1)t \Bigr)$$
be the symmetric moment curve. Then, for any $t \in \R^1$, the vectors
$$\SM_{2k}(t), \ \frac{d}{dt} \SM_{2k}(t),\ \frac{d^2}{dt^2}\SM_{2k}(t),\ 
\ldots,\ \frac{d^{2k-1}}{dt^{2k-1}} 
\SM_{2k}(t)$$
are linearly independent.
\end{lemma} 
\proof
Because of rotational invariance, 
it suffices to prove the result for $t=0$.
Let us consider the $2k$ vectors 
$$\SM_{2k}(0),\ \frac{d}{dt}\SM_{2k}(0), \ 
\ldots, \ \frac{d^{2k-1}}{d t^{2k-1}}\SM_{2k}(0),$$
that is, the vectors
\begin{eqnarray*}
a_j&=&(-1)^j \Bigl(1,\ 0,\ 3^{2j},\ 0,\ \ldots,\ 0, \ (2k-1)^{2j} \Bigr) \quad
\quad \text{and} \\
b_j&=&(-1)^j \Bigl(0,\ 1,\ 0,\ 3^{2j+1}, \ \ldots, \ (2k-1)^{2j+1},\  0 \Bigr)
\end{eqnarray*}
for $j=0, \ldots, k-1$.
It is seen then that the set of  vectors 
$\left\{a_j, b_j: j=0, \ldots, k-1\right\}$ is 
linearly independent if and only if both sets of vectors $\left\{a_j, j=0, 
\ldots, k-1\right\}$ and $\left\{b_j: j=0, \ldots, k-1\right\}$ 
are linearly independent. On the other hand, the odd-numbered 
coordinates of $(-1)^ja_j$ form the  $k \times k$
Vandermonde matrix
$$\left(\begin{matrix} 1 & 1 & 1 & \ldots & 1 \\
1 & 3^2 & 5^2 &\ldots & (2k-1)^2  \\
\ldots & \ldots & \ldots & \ldots & \ldots \\
 1 & 3^{2k-2}, & 5^{2k-2} &\ldots & (2k-1)^{2k-2} 
\end{matrix}\right) $$
while the even-numbered coordinates of $(-1)^j b_j$ 
form the $k \times k$ Vandermonde matrix

$$\left(\begin{matrix} 1 & 3 & 5 & \ldots & (2k-1) \\
1 & 3^3 & 5^3 & \ldots & (2k-1)^3 \\
\ldots & \ldots & \ldots & \ldots & \ldots \\
1 & 3^{2k-1} & 5^{2k-1} & \ldots & (2k-1)^{2k-1}
\end{matrix}\right).$$
Hence the statement follows.
\endproof

Next, we establish a curious 
property of zeros of raked trigonometric polynomials.

\begin{lemma} \label{vert}
Let
$$A(t)=c_0 + \sum_{j=1}^k a_j \sin (2j-1)t +\sum_{j=1}^k b_j \cos (2j-1)t$$
be a raked trigonometric polynomial  $A: \S^1 \longrightarrow \R$ that is 
not identically 0.
Suppose that $A$ has  $2k$ distinct roots in an arc $\Omega \subset \S^1$ 
of length less than $\pi$. Then, if $A$ has yet another root  
on $\S^1$, that root must lie in the arc 
$\Omega +\pi$.
\end{lemma}
\proof
Let us consider the derivative of $A(t)$,
$$A'(t)=\sum_{j=1}^k a_j (2j-1) \cos (2j-1)t -\sum_{j=1}^k b_j 
(2j-1) \sin (2j-1) t,$$
$A': \S^1 \longrightarrow \R$.
Substituting $z=e^{it}$, we can write 
$$A'(t)=\frac{1}{z^{2k-1}} P(z),$$
where $P(z)$ is a polynomial of degree $4k-2$, cf. Section \ref{raked-rel}.
 Hence the total number of the roots of $A'$ in $\S^1$, 
counting multiplicities, does not exceed $4k-2$.

Let $t_0,t_1\in \Omega$ be the roots of $A$ closest to 
the endpoints of $\Omega$.
By Rolle's Theorem, $A'$ has at least $2k-1$ distinct roots between 
$t_0$ and $t_1$ in $\Omega$.
Since $A'(t + \pi)=-A'(t)$, we must have another $2k-1$ distinct roots of 
$A'$ in the arc $\Omega+\pi$ between $t_0+\pi$ and $t_1+\pi$, see Figure 2.  

\begin{figure}
$$\epsffile{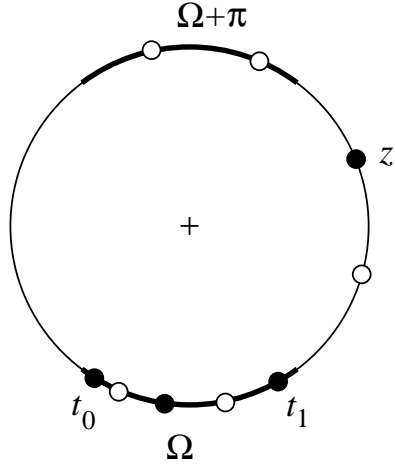}$$
\caption{The roots of $A$ (black dots) and roots of $A'$ (white dots)}
\end{figure}

Suppose that $A$ has a root $z \in \S^1$ outside of 
$\Omega \cup (\Omega+\pi)$. 
Then either $z$ lies in the open arc with the endpoints $t_0$ and 
$t_1+\pi$
or $z$ lies in the open arc with the endpoints $t_1$ and $t_0+\pi$. By Rolle's 
Theorem, $A'$ has yet another root in $\S^1$ between $t_0$ on $z$ 
in the first case, and between
$z$ and $t_1$ in the second case, which is a contradiction.
\endproof

We are now ready to prove Proposition \ref{rakegen}.

\smallskip\noindent{\it Proof of Proposition \ref{rakegen}}: 
First, we observe that for any 
$2k$ points $t_1, \ldots, t_{2k} \in \S^1$ 
there is an affine hyperplane passing through the points 
$\SM_{2k}(t_1), \ldots, \SM_{2k}(t_{2k})$ in $\R^{2k}$ 
and hence there is a non-zero raked trigonometric 
polynomial $A$ 
such that $A(t_1)=\ldots =A(t_{2k})=0$. Moreover, 
if $t_1, \ldots, t_{2k}$ are distinct and lie in an arc $\Omega$ 
of length less
than $\pi$ then the hyperplane is unique. Indeed, if the hyperplane 
is not unique then the points $\SM_{2k}(t_1), \ldots, \SM_{2k}(t_{2k})$ 
lie in an affine subspace of codimension
at least 2. Therefore, for any point 
$t_{2k+1} \in \Omega \setminus\{t_1, \ldots, t_{2k}\}$ there is 
an affine hyperplane passing through 
$\SM_{2k}(t_1), \ldots, \SM_{2k}\left(t_{2k+1}\right)$ and hence there is 
a raked polynomial that has $2k+1$ roots in $\Omega$ 
and is not identically 0, which contradicts
Lemma \ref{vert}.

Suppose now that no matter how small $\phi_k>0$ is, 
there is always an arc $\Omega \subset 
\S^1$ of length at most $\phi_k$ and a non-zero raked 
polynomial $A$ of degree $2k-1$ which has 
$2k$ distinct roots 
in $\Omega$ and at least one more root elsewhere in $\S^1$. 
By Lemma \ref{vert}, that 
remaining root must lie in the arc $\Omega + \pi$. In other words,  
for any positive integer $n$ there exists an arc $\Omega_n \subset \S^1$
of length at most $1/n$ and an affine hyperplane
$H_n$ which intersects $\SM_{2k}\left(\Omega_n\right)$ 
in $2k$ distinct points and intersects 
the set $\SM_{2k}\left(\Omega_n+\pi\right)$ as well. 
The set of all affine hyperplanes
intersecting the compact set $\SM_{2k}\left(\S^1\right)$ 
is compact in the natural topology; for example if we view
the set of affine hyperplanes in  $\R^{2k}$ as a subset of 
the Grassmannian of all (linear) 
hyperplanes in $\R^{2k+1}$. Therefore, the sequence of hyperplanes 
$H_n$ has a limit 
hyperplane $H$. By Lemma \ref{noflat}, the affine hyperplane $H$
is the $(2k-1)$th order tangent hyperplane to $\SM_{2k}\left(\S^1\right)$ 
at  some point $\SM_{2k}\left(t_0\right)$ where $t_0$ is a limit point 
of the arcs $\Omega_n$. Also, $H$ passes through the point 
$-\SM_{2k}\left(t_0\right)$. The corresponding 
trigonometric polynomial $A(t)$ is a raked polynomial 
of degree at most $2k-2$ that is not 
identically 0 and has two roots $t_0$ and 
$t_0+\pi$ with the multiplicity of $t_0$ being at least $2k$.

Let 
$$A(t)=c_0 + \sum_{j=1}^k a_j \sin (2j-1) t +\sum_{j=1}^k b_j \cos(2j-1)t.$$
Since $A(t_0)=A(t_0+\pi)=0$ we conclude that $c_0=0$. 
This, however, contradicts Lemma~\ref{noflat}
since the non-zero $2k$-vector
$$\Bigl(b_1,\  a_1,\  b_2, \ a_2, \ \ldots, b_k,\ a_k \Bigr)$$
 turns out to be orthogonal to vectors 
$$\SM_{2k}\left(t_0\right), \ \frac{d}{dt} \SM_{2k}\left(t_0\right),\  
\ldots, \ \frac{d^{2k-1}}{dt^{2k-1}} \SM_{2k}\left(t_0\right).$$
\endproof

\smallskip\noindent{\it Proof of Theorem \ref{facegen}}: 
Let $\phi_k>0$ be the number whose
existence is established in Proposition \ref{rakegen}. 
 Given $k$ distinct points $t_1, \ldots, t_k$ lying
on an arc of length at most $\phi_k$, we must present 
a raked trigonometric polynomial $A$ that has roots of multiplicity 2 
at $t_1, \ldots, t_k$ and no  other roots
on the circle. In geometric terms, 
we must present an affine hyperplane that is the first order 
tangent to the points $\SM_{2k}(t_1), \ldots, \SM_{2k}(t_k)$ 
and does not intersect $\SM_{2k}(\S^1)$ anywhere 
else. As in the proof of Proposition \ref{rakegen}, such 
a hyperplane is obtained as a limit of the affine hyperplanes which, 
for every $j=1, \ldots, k$ intersect $\SM_{2k}(\S^1)$ 
at two distinct points converging to $t_j$.
\endproof

\smallskip\noindent {\it Proof of Theorem \ref{main-thm2}: \ } 
The upper bound follows by Proposition \ref{UB2}. 
To prove the lower bound, let us consider the polytope
$\B_{2k}(X)$, where $X \subset \S^1$ is the set of $n$ 
equally spaced points ($n$ is even).
The lower bound follows by Theorem \ref{facegen}. In fact, one can show that 
$c_j(d) \geq 2^{-j-1}$: 
to obtain this inequality consider the polytope $\B_{2k}(Z)$ where 
$Z=Y \cup \left(Y+\pi\right)$ and $Y$ 
lies in an arc of length at most $\phi_k$ as defined in Theorem \ref{facegen}.
\endproof

\section{Concluding remarks}
We close the paper with two additional remarks on the face numbers
of centrally symmetric polytopes and several open questions.

\subsection{The upper half of the face vector}
Theorems \ref{main-thm1} and \ref{main-thm2} provide estimates 
on $\fmax(2k, n; j)$ --- the maximal possible number of 
$j$-faces that a cs $2k$-polytope on $n$ vertices can have --- for 
$j\leq k-1$. What can be said about $\fmax(2k, n; j)$ for $j \geq k$?
Here we prove that for every $k \leq j< 2k$, the value of $\fmax(2k, n; j)$ 
is in the order of $n^k$.
\begin{theorem}  \label{main-thm3}
Let us fix a positive even integer $d=2k$ and an integer $k \leq j < 2k$.
Then there exist $\gamma_j(d), \Gamma_j(d) >0$ such that

\begin{equation*}
\gamma_j(d) +o(1) \leq \frac{\fmax(d,n; j)}{{n \choose k}} 
 \leq \Gamma_j(d)+o(1) \quad \text{as} \quad n \longrightarrow +\infty.
\end{equation*}

\end{theorem}
\proof The upper bound estimate follows from 
the Upper Bound Theorem \cite{McMu70} 
which holds for all polytopes. 
To verify the lower bound, let us consider a cs $2k$-polytope $P_n$
on $n$ vertices that satisfies  $f_{k-1}(P_n)=\fmax(2k,n; k-1)$.
As in the proof of Proposition \ref{UB2} we can assume that $P_n$ 
is a simplicial polytope. 
Let 
$$h(P_n)=\Bigl(h_0(P_n),\  h_1(P_n),\  \ldots,\  h_{2k}(P_n)\Bigr)$$ 
be the $h$-vector of $P_n$ (see for instance 
\cite[Chapter 8]{Zieg}), that is, 
the vector whose entries are defined by the polynomial identity
$$
\sum_{i=0}^d h_i(P_n)x^{2k-i} = 
\sum_{i=0}^d f_{i-1}(P_n) (x-1)^{2k-i}.
$$
Equivalently,
\begin{equation}  \label{h-numbers}
f_{j-1}(P_n)=\sum_{i=0}^j {2k-i \choose 2k-j}h_i(P_n), 
\quad j=0, 1, \ldots, 2k.
\end{equation}
The $h$-numbers of a simplicial polytope are well-known to be
nonnegative and symmetric  \cite[Chapter 8]{Zieg}, that is,
$h_j(P_n)=h_{2k-j}(P_n)$ for $j=0,1, \ldots, 2k$.
Moreover, McMullen's proof of the UBT
implies that the $h$-numbers of any simplicial $2k$-polytope with $n$
vertices satisfy 
$$h_j \leq {n-2k+j-1 \choose j}=O(n^j), \quad \text{for} \quad 
0\leq j\leq k.
$$
Substituting these inequalities into (\ref{h-numbers}) for  
$j=k-1$ and using  that 
$$f_{k-1}(P_n)=\fmax(2k, n; k-1)=\Omega\left(n^k\right)$$
by Theorem \ref{main-thm2}, we obtain 
$$h_k(P_n)= \Omega\left(n^k\right). $$ 
Together with
nonnegativity of $h$-numbers and (\ref{h-numbers}), this implies that
$$\fmax(2k,n; j)\geq f_{j}(P_n)=\Omega\left(n^k\right) 
\quad \text{for all}  \quad k\leq j < 2k,$$
as required.
\endproof

\subsection{2-faces of $\B_6$}

We provide some additional estimates on the extent to which
$\B_6$ is 3-neighborly.

\begin{theorem}
Let $t_1, t_2, t_3\in \R$ be such that the points 
$z_1=e^{it_1}$, $z_2=e^{it_2}$, and 
$z_3=e^{it_3}$ are distinct and lie on an arc of the unit circle of
length at most $\arccos (1/8)$. 
Then the convex hull of the set $\{\SM_6(t_1),\ \SM_6(t_2),\ \SM_6(t_3)\}$
is a 2-dimensional face of $\B_6$. Consequently, 
$$\frac{\fmax(6, n; 2)}{{n \choose 3}} 
\geq 3\left(\frac{\arccos 1/8}{2\pi}\right)^2 + o(1) \approx 0.159.$$ 
\end{theorem}
\proof As in  Proposition \ref{rakegen} and Theorem~\ref{facegen}, the 
proof reduces to verifying the following statement:

\medskip
Let $z_1, \cdots, z_6 \in \C$ be distinct
points that lie on an arc of the unit circle $|z|=1$ 
of length at most $\arccos (1/8)$. Let $D(z)$ be a raked
self-inversive polynomial of degree 10  such that 
$D(z_j)=0$ for $j=1, \ldots, 6$. Then 
 none of the remaining roots of $D$ have the absolute value of 1.
\medskip 
 
Let $z_7, z_8, z_9$, and $z_{10}$ be the remaining roots of 
$D$ (some of the roots may coincide).
Let $\Phi$ be an arc of the unit circle $|z|=1$ of length
 $l\leq \arccos(1/8)$ that contains $z_1, \ldots, z_6$ and let us 
consider the line $L$ through the origin that bisects $\Phi$.
Since $D$ is a raked polynomial, we must have

\begin{equation}\label{sums}
\sum_{j=1}^{10} z_j = \sum_{j=1}^{10} z_j^3=0,
\end{equation}
 cf. Lemma \ref{roots-selfinv}.
Let  $\Sigma_1$ be  the sum of the orthogonal  projections 
of $z_1, \ldots, z_6$
onto $L$ and let $\Sigma_2$ be the sum of  the orthogonal projections of 
$z_7, \cdots, z_{10}$ onto $L$, so 
$$\Sigma_1 + \Sigma_2 =0.$$
As $\cos l \geq 1/8$, we have
$\cos (l/2)  \geq 3/4$, and hence 
\begin{equation}\label{projection}
|\Sigma_2|=|\Sigma_1|\geq 6\cdot \frac{3}{4}=\frac{9}{2}.
\end{equation}
Therefore, for at least one of the roots of $D$, say, $z_9$ 
we have $|z_9|>1$.
Then, for another root of $D$, say, $z_{10}$ we have $|z_{10}|=1/|z_9|<1$, 
cf.~Lemma \ref{roots-selfinv}.
If $|z_7|>1$ then $|z_8|<1$ and we are done.
 Hence the only remaining case to consider 
is $|z_7|=|z_8|=1$. In this case, by (\ref{projection}),
 we should have $|z_9| \geq 2$. 
Using that $z_{10}=1/\overline{z_9}$ we obtain
$$|z_9^3+z_{10}^3|=|z_9^3|+|z_{10}^3|>8=\sum_{j=1}^8 |z_j|^3 
\geq \big|\sum_{j=1}^8 z_j^3\big|,
$$
which contradicts (\ref{sums}).
\endproof

\subsection{Open questions}

There are several natural questions that we have not be able to answer so far.
\begin{itemize}
\item 
It seems plausible that $\psi_k$ in Theorem \ref{lower1} satisfies
$$\psi_k=\frac{2k-2}{2k-1} \pi, $$
 but we are unable to prove that.
\item It also seems plausible that in Theorem \ref{main-thm2} 
for any fixed $j$ we have
$$\lim_{d \longrightarrow +\infty} c_j(d)=1,$$
but we are also unable to prove that.
\item
We do not know what is the best value of 
$\phi_k$ in Theorem \ref{facegen} for $k>2$ nor
the values of $c_j(d)$ in Theorem \ref{main-thm2}.
\item
The most intriguing question is, of course, 
whether the class of polytopes $\B_{2k}(X)$
indeed provides (asymptotically or even exactly) 
polytopes with the largest number of faces
among all centrally symmetric polytopes with a given number of vertices. 
\end{itemize}

\section*{Acknowledgments} 
The authors are grateful to J.E. Goodman, R. Pollack, 
and J. Pach, the organizers of the 
 AMS - IMS - SIAM Summer Research Conference ``Discrete and 
Computational Geometry - twenty years later'' 
(Snowbird, June 2006), where this project 
started, and to L. Billera for encouragement.


\begin{thebibliography}{999}
\bibitem{Barv} A.~Barvinok, {\em A Course in Convexity}, 
Graduate Studies in Mathematics, {\bf 54}, American Mathematical Society, 
Providence, RI, 2002.

\bibitem {BiLee}
L.~J.~Billera and C.~W.~Lee,
``A proof of the sufficiency of McMullen's conditions
for $f$-vectors of simplicial convex polytopes",
J.~Comb.~Theory Ser.~A {\bf 31} (1981), 237--255.  

\bibitem{Car} C.~ Caratheodory, ``\"Uber den Variabilitatsbereich det 
Fourierschen Konstanten von Positiven harmonischen Furktionen", 
Ren.~Circ.~Mat.~Palermo {\bf 32} (1911), 193-217. 

\bibitem{DanzGr} L.~Danzer and B.~Gr\"unbaum, 
``\"Uber zwei Probleme bez\"uglich konvexer K\"orper von P.~Erd\"os 
und von V.~L.~Klee'' 
(German), Math.~Z. {\bf 79} (1962), 95--99.

\bibitem{Dono1} D.~L.~Donoho, 
``High-dimensional centrosymmetric polytopes with neighborliness 
proportional to dimension", Discrete Comput.~Geom. {\bf 35} (2006),  
617--652.

\bibitem{Dono2} D.~L.~Donoho, ``Neighborly polytopes and sparse
solutions of underdetermined linear equations", preprint
(2004).

\bibitem{DonoTan} D.~L.~Donoho and J.~Tanner, 
``Counting faces of randomly-projected polytopes 
when the projection radically lowers dimension", preprint, math.MG/0607364.  

\bibitem{Gale} D.~Gale,
``Neighborly and cyclic polytopes", Proc.~Sympos.~Pure Math., Vol. VII,
Amer.~Math.~Soc., Providence, R.I., 1963,  pp.~225--232. 

\bibitem{Grunb} B.~Gr\"unbaum, {\em Convex polytopes}, second edition 
(prepared and with a preface by V.~Kaibel, V.~Klee and G.~M.~Ziegler),
 Graduate Texts in Mathematics, {\bf 221}, Springer-Verlag, New York, 2003.

\bibitem{GrMotz} B.~Gr\"unbaum and T.~S.~Motzkin, ``On polyhedral gaps",
Proc.~Sympos.~Pure Math., Vol. VII,
Amer.~Math.~Soc., Providence, R.I., 1963,  pp.~285-290.

\bibitem{Ka} G.~O.~H.~Katona, ``A theorem of finite sets",
Theory of graphs (Proc.~Colloq., Tihany, 1966), Academic Press, 
New York, 1968, pp. 187--207.

\bibitem{Kr} J.~B.~Kruskal, ``The number of simplices in a complex",  
Mathematical optimization techniques,  Univ.~of California Press, Berkeley, 
Calif., 1963, pp. 251--278. 

\bibitem{LinNov} N.~Linial and I.~Novik, ``How neighborly can a centrally 
symmetric polytope be?",  Discrete Comput.~Geometry {\bf 36} (2006),
273--281.

\bibitem{McMu70}
P.~McMullen, ``The maximum numbers of faces of a convex polytope",
  Mathematika {\bf 17} (1970), 179--184.

\bibitem{Motz} T.~S.~Motzkin, ``Comonotone curves and polyhedra".
Bull.~Amer.~Math.~Soc. {\bf 63} (1957), 35.

\bibitem{RudelVersh} M.~Rudelson and R.~Vershynin,
``Geometric approach to error correcting codes and reconstruction of 
signals", Int.~Math.~Res.~Not. {\bf 64} (2005), 4019--4041.

\bibitem{Small} T.~Sheil-Small, {\em Complex polynomials},
 Cambridge Studies in Advanced Mathematics, {\bf 75}, 
Cambridge University Press, Cambridge, 2002. 

\bibitem{Smil85} Z.~Smilansky, ``Convex hulls of generalized moment curves", 
Israel J.~Math. {\bf 52} (1985), 115--128.

\bibitem{Smil90} Z.~Smilansky, ``Bi-cyclic $4$-polytopes",
Israel J.~Math. {\bf 70} (1990),  82--92.

\bibitem{St80}
R.~Stanley, ``The number of faces of simplicial convex polytopes",
Adv.~Math. {\bf 35} (1980), 236--238.

\bibitem{Zieg}  
G.~M.~Ziegler,  
{\em Lectures on Polytopes},
 Graduate Texts in Mathematics, {\bf 152}, Springer-Verlag, New York, 1995. 

\end{thebibliography}
\end{document}